\documentclass[11pt,leqno]{amsart}
\usepackage{amssymb,verbatim,enumerate,ifthen,hyperref}
\usepackage[mathscr]{eucal}
\usepackage[utf8]{inputenc}
\usepackage[T1]{fontenc}
\def\N{\mathbb{N}}
\def\R{\mathbb{R}}

\def\M{\mathscr{M}}
\def\A{\mathscr{A}}
\def\H{\mathscr{H}}
\def\E{{\mathscr{E}}}

\def\sign{\mathop{\mbox{\rm sign}}\nolimits}
\long\def\comment#1{}

\newtheorem{theorem}{Theorem}[section]
\newtheorem*{theorem*}{Theorem}
\def\Thm#1#2{\ifthenelse{\equal{#1}{*}}{\begin{theorem*}#2\end{theorem*}}
             {\begin{theorem}\label{T#1}#2\end{theorem}}}
\newtheorem{Atheorem}{Theorem}

\def\thm#1{Theorem~\ref{T#1}}
\newtheorem{proposition}[theorem]{Proposition}
\newtheorem*{proposition*}{Proposition}
\def\Prp#1#2{\ifthenelse{\equal{#1}{*}}{\begin{proposition*}#2\end{proposition*}}
{\begin{proposition}\label{P#1}#2\end{proposition}}}
\def\prp#1{Proposition~\ref{P#1}}

\newtheorem{corollary}[theorem]{Corollary}
\newtheorem*{corollary*}{Corollary}
\def\Cor#1#2{\ifthenelse{\equal{#1}{*}}{\begin{corollary*}#2\end{corollary*}}
             {\begin{corollary}\label{C#1}#2\end{corollary}}}
\def\cor#1{Corollary~\ref{C#1}}

\newtheorem{lemma}[theorem]{Lemma}
\newtheorem*{lemma*}{Lemma}
\def\Lem#1#2{\ifthenelse{\equal{#1}{*}}{\begin{lemma*}#2\end{lemma*}}
             {\begin{lemma}\label{L#1}#2\end{lemma}}}
\def\lem#1{Lemma~\ref{L#1}}

\theoremstyle{definition}
\newtheorem{remark}[theorem]{Remark}
\newtheorem*{remark*}{Remark}
\def\Rem#1#2{\ifthenelse{\equal{#1}{*}}{\begin{remark}\rm #2\end{remark}}
             {\begin{remark}\label{R#1}\rm #2\end{remark}}}

\newtheorem{example}[theorem]{Example}
\newtheorem*{example*}{Example}
\def\Exa#1#2{\ifthenelse{\equal{#1}{*}}{\begin{example*}\rm #2\end{example*}}
             {\begin{example}\label{Ex#1}\rm #2\end{example}}}

\def\eq#1{{\rm(\ref{E#1})}}
\def\Eq#1#2{\ifthenelse{\equal{#1}{*}}
  {\begin{equation*}\begin{aligned}#2\end{aligned}\end{equation*}}
  {\begin{equation}\begin{aligned}\label{E#1}#2\end{aligned}\end{equation}}}

\begin{document}
\begin{flushright}
%\textit{Submitted to:} 
\end{flushright}
\vspace{5mm}

\date{\today}

\title[Approximately monotone and approximately H\"older functions]{Characterization of approximately monotone \\ and approximately H\"older functions}

\author[A. R. Goswami]{Angshuman R. Goswami}
\address[A. R. Goswami]{Doctoral School of Mathematical and Computational Sciences, University of Debrecen, 
H-4002 Debrecen, Pf.\ 400, Hungary}
\author[Zs. P\'ales]{Zsolt P\'ales}
\address[Zs. P\'ales]{Institute of Mathematics, University of Debrecen, 
H-4002 Debrecen, Pf.\ 400, Hungary}
\email{\{pales,angshu\}@science.unideb.hu}


\thanks{The research of the second author was supported by the K-134191 NKFIH Grant and the 2019-2.1.11-TÉT-2019-00049, the EFOP-3.6.1-16-2016-00022 and the EFOP-3.6.2-16-2017-00015 projects. The last two projects are co-financed by the European Union and the European Social
Fund.
}

\begin{abstract}
A real valued function $f$ defined on a real open interval $I$ is called $\Phi$-monotone if, for all $x,y\in I$ with $x\leq y$ it satisfies
$$
  f(x)\leq f(y)+\Phi(y-x),
$$
where $\Phi:[0,\ell(I)[\,\to\mathbb{R}_+$ is a given nonnegative error function, where $\ell(I)$ denotes the length of the interval $I$. If $f$ and $-f$ are simultaneously $\Phi$-monotone, then $f$ is said to be a $\Phi$-H\"older function.
In the main results of the paper, using the notions of upper and lower interpolations, we establish a characterization for both classes of functions. This allows one to construct $\Phi$-monotone and $\Phi$-H\"older functions from elementary ones, which could be termed the building blocks for those classes. In the second part, we deduce Ostrowski- and Hermite--Hadamard-type inequalities from the $\Phi$-monotonicity and $\Phi$-H\"older properties, and then we verify the sharpness of these implications. We also establish implications in the reversed direction.
\end{abstract}

\maketitle

\section{INTRODUCTION}

The notions of approximately monotone, submonotone and approximately H\"older functions were coined up in several papers (cf.\ \cite{Csa60,DanGeo04,ElsPeaRob96,MakPal12a,NgaPen07,Pal03a}) and have applications in nonsmooth and convex analysis and optimization theory, and also in the theory of functional equations and inequalities. Motivated by the applicability of such concepts, in our former paper \cite{GosPal20}, we introduced the concepts of approximately monotone and approximately H\"older functions and established their basic properties. We now recall the terminologies and notations introduced in \cite{GosPal20}.

Let $I$ be a nonempty open real interval throughout this paper and let $\ell(I)\in\,]0,\infty]$ denote its length. The symbols $\R$ and $\R_+$ denote the sets of real and nonnegative real numbers, respectively. The class of all functions $\Phi:[0,\ell(I)[\,\to\R_+$, called error functions, will be denoted by $\E(I)$. We define two concepts related to an error function $\Phi\in\E(I)$.

A function $f:I\to\R$ is called \emph{$\Phi$-monotone} if, for all $x,y\in I$ with $x\leq y$, 
\Eq{H1}{
  f(x)\leq f(y)+\Phi(y-x).
}
If this inequality is satisfied with the identically zero error function $\Phi$, then we say that $f$ is \emph{monotone (increasing)}. The class of all $\Phi$-monotone functions on $I$ will be denoted by $\M_\Phi(I)$.

A function $f:I\to\R$ is called \emph{$\Phi$-H\"older} if, for all $x,y\in I$, 
\Eq{H2}{
  |f(x)-f(y)|\leq\Phi(|x-y|).
}
The class of all $\Phi$-H\"older functions on $I$ will be denoted by $\H_\Phi(I)$.

In \cite{GosPal20}, we showed that the classes $\M_\Phi(I)$ and $\H_\Phi(I)$ are convex and closed with respect to the pointwise supremum, i.e., if $\{f_\gamma:I\to\R\mid\gamma\in\Gamma\}$ is a subfamily of $\M_\Phi(I)$ [resp.\ $\H_\Phi(I)$] with a pointwise supremum $f:I\to\R$, i.e., 
\Eq{*}{
  f(x)=\sup_{\gamma\in\Gamma} f_\gamma(x)\qquad(x\in I),
}
then $f\in\M_\Phi(I)$ [resp.\ $\H_\Phi(I)$]. Analogously, $\M_\Phi(I)$ and $\H_\Phi(I)$ are also closed with respect to the pointwise infimum. Additionally, the class $\H_\Phi(I)$ is centrally symmetric, i.e., $\H_\Phi(I)$ is closed with respect to multiplication by $(-1)$. One can also observe that
\Eq{*}{
\H_\Phi(I)=\M_\Phi(I)\cap(-\M_\Phi(I)).
}
In \cite{GosPal20}, we proved that an optimal error functions for $\Phi$-monotonicity and $\Phi$-H\"older property must be subadditive and absolutely subadditive, respectively. We recall, that an error function $\Phi\in\E(I)$ is \emph{subadditive} if, for all $u,v\in\R_+$ with $u+v<\ell(I)$, the inequality
\Eq{*}{
  \Phi(u+v)\leq \Phi(u)+\Phi(v)
}
holds. Similarly, $\Phi$ is termed as \emph{absolutely subadditive} if, for all $u,v\in\R$ with $|u|,|v|,|u+v|<\ell(I)$, the inequality
\Eq{*}{
  \Phi(|u+v|)\leq \Phi(|u|)+\Phi(|v|)
}
is satisfied. It is clear that absolutely subadditive functions are automatically subadditive. On the other hand, as we have proved it in \cite{GosPal20}, increasingness and subadditivity imply absolute subadditivity. 
In \cite{GosPal20}, we also established a formula for the lower and for the upper $\Phi$-monotone and $\Phi$-H\"older envelopes. Furthermore, we introduced a generalization of the classical notion of total variation and proved an extension of the Jordan Decomposition Theorem known for functions of bounded total variations.

The aim of this paper is twofold. In the first part, using the notions of upper and lower interpolations, we establish a characterization for both classes of functions. This allows one to construct $\Phi$-monotone and $\Phi$-H\"older functions from elementary ones, which could be termed as the building blocks for those classes. In addition to these, we describe the solution to a two variable triangle inequality-type functional equation which appears in the characterization of $\Phi$-monotonicity. In the second part, we deduce Ostrowski- and Hermite--Hadamard-type inequalities from the $\Phi$-monotonicity and $\Phi$-H\"older properties, and then we verify the sharpness of these implications. We also establish implications in the reversed direction.

\section{Characterization of $\Phi$-monotone functions}

In what follows, we construct a large class of elementary $\Phi$-monotone functions provided that $\Phi$ is a subadditive and nondecreasing error function whose members will turn out to be the building blocks of $\Phi$-monotone functions. It turns out (cf. \cite[Proposition 3.2]{GosPal20}), that if $\Phi$ is of the form $\Phi(t)=ct^p$, where $c\geq0$ and $p\in\,]0,1]$, then $\Phi$ is a subadditive and nondecreasing error function.

If $\Phi\in\E(I)$, $h:I\to[-\infty,\infty]$ and $p\in I$, then define the functions $h_p,h^p:I\to[-\infty,\infty]$ as follows:
\Eq{hp}{
  h_p(x):=\begin{cases}
         h(x) &\mbox{if } x\leq p,\\%[2mm]
         h(p)-\Phi(x-p) &\mbox{if } p<x,
        \end{cases}
    \qquad\mbox{and}\qquad
 h^p(x):=\begin{cases}
         h(p)+\Phi(p-x) &\mbox{if } x<p,\\%[2mm]
         h(x) &\mbox{if } p\leq x.
        \end{cases}
}

\Prp{1}{Let $\Phi\in\E(I)$ be subadditive and nondecreasing and $h:I\to[-\infty,\infty]$ be nondecreasing. Then, for all $p\in I$, the functions $h_p$ and $h^p$ defined by \eq{hp} are $\Phi$-monotone.}

\begin{proof} Let $p\in I$ be fixed. To show that $h_p$ is $\Phi$-monotone, let $x,y\in I$ with $x<y$. We now distinguish three subcases according to the position of $p$ with respect to $x$ and $y$.
\\
If $x<y\leq p$, then
\Eq{*}{
  h_p(x)=h(x)\leq h(y)=h_p(y)\leq h_p(y)+\Phi(y-x).
}
If $x\leq p<y$, then the increasingness of $\Phi$ implies
\Eq{*}{
h_p(x)=h(x)\leq h(p)= h_p(y)+\Phi(y-p)\leq h_p(y)+\Phi(y-x).
}
Finally, if $p<x<y$, then the subadditivity of $\Phi$ results
\Eq{*}{
h_p(x)=h(p)-\Phi(x-p)=h_p(y)+\Phi(y-p)-\Phi(x-p)\leq h_p(y)+\Phi(y-x).
}
The proof of the $\Phi$-monotonicity of $h^p$ is analogous, therefore omitted.
\end{proof}

If $\Phi\in\E(I)$, $f:I\to\R$ and $p\in I$, then we say that \emph{$f$ can be interpolated at $p$ by a $\Phi$-monotone function from below [resp. from above]} if there exits a $\Phi$-monotone function $h:I\to\R$ such that $h(p)=f(p)$ and $h\leq f$ [resp. $f\leq h$].

\Prp{A}{Let $\Phi\in\E(I)$ be a subadditive and nondecreasing function, let $f:I\to\R$ and $p\in I$ be fixed. Then $f$ can be interpolated at $p$ by a $\Phi$-monotone function from below if and only if, for all $x\in I$,
 \Eq{cond}{
    -\infty<\inf_{[x,p]}f \quad\mbox{if }x\leq p \qquad\mbox{and}\qquad 
    f(p)\leq f(x)+\Phi(x-p) \quad\mbox{if }p<x.
 }
Analogously, $f$ can be interpolated at $p$ by a $\Phi$-monotone function from above if and only if, for all $x\in I$,
 \Eq{cond2}{
    f(x)\leq f(p)+\Phi(p-x) \quad\mbox{if }x<p \qquad\mbox{and}\qquad 
    \sup_{[p,x]}f<+\infty \quad\mbox{if }p\leq x.
 }
}

\begin{proof} Suppose first that there exists a $\Phi$-monotone function $h:I\to\R$ such that $h(p)=f(p)$ and $h\leq f$. Let $x\in I$ be arbitrary. 

If $x\leq p$ and $t\in[x,p]$, then by the $\Phi$-monotonicity of $h$ and the nondecreasingness of $\Phi$, we get
\Eq{*}{
h(x)\leq h(t)+\Phi(t-x)\leq h(t)+\Phi(p-x) \leq f(t)+\Phi(p-x),
}
which implies that 
\Eq{*}{
  \inf_{[x,p]}f\geq h(x)-\Phi(p-x)>-\infty.
}
If $p<x$, then $f(p)=h(p)\leq h(x)+\Phi(x-p)\leq f(x)+\Phi(x-p)$ proving the second part of condition \eq{cond}.

To prove the sufficiency, assume that condition \eq{cond} holds and define the function $h:I\to[-\infty,\infty]$ by
\Eq{*}{
  h(x):=\begin{cases}
         \displaystyle{\inf_{[x,p]}f} &\mbox{if } x<p,\\[4mm]
         f(p) &\mbox{if } x=p,\\[1mm]
         \displaystyle{\sup_{[p,x]}f} &\mbox{if } p<x.
        \end{cases}
}
Then one can easily check that $h$ is nondecreasing, $h(p)=f(p)$ and, by the first part of \eq{cond}, $h(x)$ is finite if $x\leq p$. Therefore, the function $h_p$ has finite values everywhere. According to \prp{1}, $h_p$ is $\Phi$-monotone and $h_p(p)=h(p)=f(p)$. Thus, it remains to show that $h_p\leq f$. Indeed, if $x\leq p$,
then $h_p(x)=\inf_{[x,p]}f\leq f(x)$. If $p<x$, then, by the second part of \eq{cond}, 
$h_p(x)=h(p)-\Phi(x-p)=f(p)-\Phi(x-p)\leq f(x)$ and the proof is completed.

The second part of the proposition can be verified in an analogous manner.
\end{proof}

\Thm{1}{Let $\Phi\in\E(I)$ be a subadditive and nondecreasing function and $f:I\to\R$. Then the following assertions are equivalent.
\begin{enumerate}[(i)]
 \item $f$ is $\Phi$-monotone.
 \item There exists a function $H:I\times I\to\R$ such that $H$ satisfies the functional equations 
 \Eq{FEH}{
   \min(H(x,y),H(y,z))=H(x,z) \qquad\mbox{and}\qquad \max(H(z,y),H(y,x))=H(z,x)
 }
 for all $x,y,z\in I$ with $x\leq y\leq z$ and, for all $p\in I$, the function $h:=H(\cdot,p)$ is nondecreasing and $h_p$ and $h^p$ are $\Phi$-monotone interpolations for $f$ at $p$ from below and from above, respectively.
 \item For every $p\in I$, there exists a nondecreasing function $h:I\to\R$ such that $h_p$ is a $\Phi$-monotone interpolation of $f$ at $p$ from below.
 \item For every $p\in I$, there exists a nondecreasing function $h:I\to\R$ such that $h^p$ is a $\Phi$-monotone interpolation of $f$ at $p$ from  above.
\end{enumerate}}

\begin{proof}
	{
		$(i)\Rightarrow(ii)$: Assume that $f$ is $\Phi$-monotone and let's define the function  $H:I\times I\to\R$ by
		\Eq{*}
		{H(x,y)=\begin{cases}
		        \displaystyle{\inf_{[x,y]}}f &\mbox{if } x<y, \\[4mm]
                f(x) &\mbox{if } x=y, \\[1mm]
               \displaystyle{ \sup_{[y,x]}}f &\mbox{if } x>y.
		        \end{cases}
		}
		By its $\Phi$-monotonicity, the function $f$ can be interpolated from below and from above by a $\Phi$-monotone function, thus, \prp{A} implies that the values of $H$ are finite.
		
		Then, for any $x,y,z \in I$ with $x\leq y\leq z$, we have
		\Eq{*}{
			\min(H(x,y),H(y,z))=\min\bigg(\inf_{[x,y]}f,\inf_{[y,z]}f\bigg)=\inf_{[x,z]}f=H(x,z).
        }
		Analogously, we can also establish the second equality:
		\Eq{*}{
			\max(H(z,y),H(y,x))=\max\bigg(\sup_{[y,z]}f,\sup_{[x,y]}f\bigg)=\sup_{[x,z]}f=H(z,x). 
        }
        Thus, we have shown that $H$ satisfies the functional equations of assertion (ii).
		
		Let $p \in I$ be fixed. To show the nondecreasingness of $h:=H(\cdot,p)$, assume $x,y\in I$ with $x<y$. If $x<y\leq p$, 
		Then, by the defination of $h$, we have
		\Eq{*}
		{h(x)=H(x,p)=\inf_{[x,p]}f\leq \inf_{[y,p]}f=H(y,p)=h(y).
		}
		If $x\leq p<y$, 
		\Eq{*}
		{h(x)=H(x,p)=\inf_{[x,p]}f\leq f(p)\leq \sup_{[p,y]}f=H(y,p)=h(y).
		}
		Finally, if $p<x<y$, then
		\Eq{*}
		{h(x)=H(x,p)=\sup_{[p,x]}f\leq \sup_{[p,y]}f=H(y,p)=h(y).
		}
		This completes the proof of the monotonicity.
		
		Finally, we just need to show that $h_p$ is a $\Phi$-monotone interpolation of $f$ at $p$ from below.
		\Eq{*}{
			h_p\leq f\leq h^p \qquad {and} \qquad h_p(p)=f(p)=h^p(p).
		}
		For the first inequality, let $x\in I$ be arbitrary. If $x\leq p$, then 
		\Eq{*}
		{h_p(x)=h(x)=H(x,p)=\inf_{[x,p]}f\leq f(x).
		}
		If $p<x$, then the $\Phi$-monotonicity of $f$ yields
		\Eq{*}
		{h_p(x)=h(p)-\Phi(x-p)=H(p,p)-\Phi(x-p)
		=f(p)-\Phi(x-p)\leq f(x).
		}
		Therefore, $h_p\leq f$ holds on $I$. The equality $h_p(p)=f(p)$ is obvious. 
		
		Analogously, one can see that $h^p$ is a $\Phi$-monotone interpolation of $f$ at $p$ from above.
		
		$(ii)\Rightarrow(iii),(iv):$ Assume that $H:I\times I\to\R$ satisfies the conditions of assertion (ii). For any $p\in I$, define $h:=H(\cdot,p)$. Then $h_p$ and $h^p$ are 
		$\Phi$-monotone interpolation of $f$ at $p$ from below and above, respectively. Hence assertions $(iii)$ and $(iv)$ are valid.
		
        The proof of the last two implications is based on the result of our previous paper \cite{GosPal20} which states that the family of $\Phi$-monotone functions is closed under pointwise supremum and infimum.
        
        $(iii)\Rightarrow (i)$ If we assume that $f$ admits a $\Phi$-monotone interpolation from below at every point $p\in I$, then $f$ is the pointwise supremum of  a family of such $\Phi$-monotone functions, and therefore, $f$ itself is a $\Phi$-monotone function.
        
        $(iv)\Rightarrow (i)$ If we assume that $f$ admits a $\Phi$-monotone interpolation from above at every point $p\in I$, then $f$ is the pointwise infimum of a family of such $\Phi$-monotone functions, and therefore, $f$ itself is a $\Phi$-monotone function.}
\end{proof}

\section{The functional equation \eq{FEH}}

In what follows we will investigate the functional equation \eq{FEH} more closely.

\Thm{FEH}{Let $H:I\times I\to\R$ and assume that $H$ satisfies \eq{FEH} for all $x,y,z\in I$ with $x\leq y\leq z$. Then
\Eq{H=}{
  H(x,y)\leq \inf_{t\in[x,y]} H(t,t) \qquad\mbox{and}\qquad
  H(y,x)\geq \sup_{t\in[x,y]} H(t,t)
}
for all $x,y\in I$ with $x\leq y$. If, in addition, $H$ is continuous at the diagonal $\{(t,t)\mid t\in I$, then the above inequalities can be replaced by equality.}

\begin{proof} 
First we prove that, for all $n\in\N$ and elements $t_0\leq t_1\leq\dots\leq t_{n-1}\leq t_n$ in $I$, the function $H$ satisfies the following generalization of \eq{FEH}:
\Eq{FEH+}{
  \min(H(t_0,t_1),\dots,H(t_{n-1},t_n))&=H(t_0,t_n)
  \qquad\mbox{and}\qquad \\
  \max(H(t_n,t_{n-1}),\dots,H(t_1,t_0))&=H(t_n,t_0).
}
This statement is trivial if $n=1$. We prove the statement for $n>1$ by induction. For $n=2$, the equalities in \eq{FEH+} are equivalent to \eq{FEH}. 
Assume now that \eq{FEH+} holds for some $n\geq2$ and let $t_0\leq t_1\leq\dots\leq t_{n-1}\leq t_n\leq t_{n+1}$ in $I$. Then, applying the inductive hypothesis and then \eq{FEH}, we get
\Eq{*}{
  \min(H(t_0,t_1),\dots,H(t_{n-1},t_n),H(t_n,t_{n+1})
  &=\min\big(\min(H(t_0,t_1),\dots,H(t_{n-1},t_n)),H(t_n,t_{n+1})\big)\\
  &=\min(H(t_0,t_n),H(t_n,t_{n+1}))=H(t_0,t_{n+1}).
}
This proves the first equality in \eq{FEH+} for the case $n+1$. The proof of the second equality is completely analogous.

To prove the first inequality in \eq{H=}, let $x\leq y$ and apply \eq{FEH+} for 
$(t_0,t_1,t_2,t_3)=(x,t,t,y)$ and for $(t_0,t_1,t_2)=(x,t,y)$ where $t\in[x,y]$ is arbitrary. We obtain
\Eq{*}{
  H(x,y)&=\min(H(x,t),H(t,t),H(t,y))\\
  &=\min(\min(H(x,t),H(t,y)),H(t,t))
  =\min(H(x,y),H(t,t))
}
This implies that $H(x,y)\leq H(t,t)$ for all $t\in[x,y]$ and completes the proof of the first inequality in \eq{H=}. The proof of the second inequality is similar.

Now assume that $H$ is continuous at the diagonal $\{(t,t)\mid t\in I$.
Let $x,y\in I$ with $x\leq y$. If $x=y$, then $[x,y]=\{x\}$ and bothe equalities in \eq{FEH} hold. Assume now that $x<y$ and let $x_0:=x$ and $y_0:=y$, and suppose that we have constructed $x_n<y_n$ such that 
\Eq{Hn}{
  [x_n,y_n]\subseteq [x_{n-1},y_{n-1}],\qquad y_n-x_n=2^{-n}(y-x)
  \qquad\mbox{and}\qquad H(x_n,y_n)= H(x_{n-1},y_{n-1}).
}
Using \eq{FEH}, we now have that
\Eq{*}{
  \min(H(x_n,\tfrac{x_n+y_n}{2}),H(\tfrac{x_n+y_n}{2},y_n))=H(x_n,y_n).
}
If $H(x_n,\tfrac{x_n+y_n}{2})=H(x_n,y_n)$ then define $x_{n+1}:=x_n$ and $y_{n+1}:=\tfrac{x_n+y_n}{2}$, otherwise let $x_{n+1}:=\tfrac{x_n+y_n}{2}$ and $y_{n+1}:=y_n$.
Then \eq{Hn} holds for $n+1$ instead of $n$.

According to the Cantor Intersection Theorem, the intersection of the sequence of intervals $[x_n,y_n]$ equals a singleton $\{t\}$, where $t\in[x,y]$ and $t$ is the common limit of the sequences $(x_n)$ and $(y_n)$. By the continuity of $H$ at $(t,t)$, it follows that the sequence $(H(x_n,y_n))$ converges to $H(t,t)$. Therefore, $H(t,t)=H(x_0,y_0)=H(x,y)$, which establishes the first equality in \eq{H=}. The proof of the second equality is similar, therefore, it is omitted.
\end{proof}

We note that the last assertion of the theorem is not valid without continuity. Indeed, define $H:I\times I\to\R$ by $H(x,y)=\sign(x-y)$. Then $H$ increases in its first variable, decreases in the second variable, and it is easy to see that $H$ satisfies the two equalities in \eq{FEH} for all $x\leq y\leq z$. On the other hand, \eq{H=} holds with strict inequalities.

It is also worth observing that a function $H:I\times I\to\R$ satisfies \eq{FEH} if and only if $G:=-H$ fulfils 
 \Eq{FEG}{
   \max(G(x,y),G(y,z))=G(x,z) \qquad\mbox{and}\qquad \min(G(z,y),G(y,x))=G(z,x)
 }
for all $x\leq y\leq z$ in $I$. Now \thm{FEH}, with the transformation $H:=-G$, implies the following assertion.

\Cor{FEH2}{Let $G:I\times I\to\R$ and assume that $G$ satisfies \eq{FEG} for all $x,y,z\in I$ with $x\leq y\leq z$. Then
\Eq{G=}{
  G(x,y)\geq \sup_{t\in[x,y]} G(t,t) \qquad\mbox{and}\qquad
  G(y,x)\leq \inf_{t\in[x,y]} G(t,t)
}
for all $x,y\in I$ with $x\leq y$. If, in addition, $G$ is continuous at the diagonal $\{(t,t)\mid t\in I\}$, then the above inequalities can be replaced by equality.}

\section{Characterization of $\Phi$-H\"older functions}

If $\Phi\in\E(I)$, $f:I\to\R$ and $p\in I$, then we say that \emph{$f$ can be interpolated at $p$ by a $\Phi$-H\"older function from below [resp. from above]} if there exits a $\Phi$-H\"older function $h:I\to\R$ such that $h(p)=f(p)$ and $h\leq f$ [resp. $f\leq h$].

In what follows, given an error function $\Phi\in\E(I)$ and $p\in I$, we define the function $\Phi_p:I\to\R$ by
\Eq{Pp}{
  \Phi_p(x):=\Phi(|x-p|) \qquad(x\in I).
}

\Prp{HH}{Let $\Phi\in\E(I)$ be an absolutely subadditive function with $\Phi(0)=0$. Then, for all $p\in I$, the function $\Phi_p$ is $\Phi$-H\"older on $I$.}

\begin{proof} Let $p\in I$ and let $x,y\in I$. Then, using the absolute subadditivity of $\Phi$, it immediately follows that
\Eq{*}{
  |\Phi_p(x)-\Phi_p(y)|=|\Phi(|x-p|)-\Phi(|y-p|)|\leq\Phi(|y-x|).
}
This shows that $\Phi_p$ is $\Phi$-H\"older on $I$.
\end{proof}

\Prp{AH}{Let $\Phi\in\E(I)$ be an absolutely subadditive function with $\Phi(0)=0$, let $f:I\to\R$ and $p\in I$ be fixed. Then $f$ can be interpolated at $p$ by a $\Phi$-H\"older function from below if and only if, for all $x\in I$,
 \Eq{condH}{
    f(p)\leq f(x)+\Phi(|x-p|).
 }
Analogously, $f$ can be interpolated at $p$ by a $\Phi$-H\"older function from above if and only if, for all $x\in I$,
 \Eq{cond2H}{
    f(x)\leq f(p)+\Phi(|x-p|).
 }
}

\begin{proof} Suppose first that there exists a $\Phi$-H\"older function $h:I\to\R$ such that $h(p)=f(p)$ and $h\leq f$. Let $x\in I$ be arbitrary. Then
\Eq{*}{
f(p)=h(p)\leq h(x)+\Phi(|x-p|)\leq f(x)+\Phi(|x-p|)
}
verifying condition \eq{condH}.

To prove the sufficiency, assume that \eq{condH} holds for all $x\in I$ and define $h:=f(p)-\Phi_p$. According to \prp{HH}, it follows that $h$ is $\Phi$-H\"older. On the other hand, the condition $\Phi(0)=0$ implies that $h(p)=f(p)$.

It remains to show that $h\leq f$. Indeed, if $x\in I$, then \eq{condH} implies 
$h(x)=f(p)-\Phi(|x-p|)\leq f(x)$ and the proof is completed.

The verification of the second assertion is analogous.
\end{proof}

\Thm{BH}{Let $\Phi\in\E(I)$ be absolutely subadditive function with $\Phi(0)=0$ and let $f:I\to\R$. Then the following assertions are equivalent to each other:
\begin{enumerate}[(i)]
 \item $f$ is $\Phi$-H\"older.
 \item For every $p\in I$, the functions $f(p)-\Phi_p$ and $f(p)+\Phi_p$ are $\Phi$-H\"older interpolations of $f$ at $p$ from below and above, respectively.
 \item For every $p\in I$, $f$ possesses a $\Phi$-H\"older interpolation from below. 
 \item For every $p\in I$, $f$ possesses a $\Phi$-H\"older interpolation from above.
\end{enumerate}
}

\begin{proof}
$(i)\Rightarrow (ii)$ Let $f$ be a $\Phi$-H\"older function. Then, trivially, $f$ possesses $\Phi$-H\"older interpolation at $p$ from below and from above. 
Thus, by \prp{AH}, the inequalities \eq{condH} and \eq{cond2H} are valid. Therefore we have that $f(p)-\Phi_p\leq f\leq f(p)+\Phi_p$. On the other hand, these two functions are $\Phi$-H\"older and interpolate $f$ at $p$. This proves $(ii)$.  

The implication $(ii)\Rightarrow (iii)$ and $(ii)\Rightarrow (iv)$ are obvious.

The proof of the implications $(iii)\Rightarrow (i)$ and $(iv)\Rightarrow (i)$ is based on a result of our previous paper \cite{GosPal20} which states that the family of $\Phi$-H\"older functions is closed under pointwise supremum and infimum. If we assume that $f$ admits a $\Phi$-H\"older interpolation from below (resp.\ from above) at every point $p\in I$, then $f$ is the pointwise supremum (resp.\ infimum) of a family of $\Phi$-H\"older functions, and therefore, $f$ itself is a $\Phi$-H\"older function.
\end{proof}

\section{Hermite--Hadamard-type inequalities for $\Phi$-monotone functions}

In the sequel, a function defined on an interval will be called locally integrable if it has a finite Lebesgue integral over every compact subinterval of its domain. 

For the description of our subsequent results, we now introduce the following notation and terminology: If $a,b\in I$ then the convex hull of $\{a,b\}$, i.e., the smallest interval containing $a$ and $b$, will denoted by $\langle a,b\rangle$. If, additionally, $f:\langle a,b\rangle\to\R$ is Lebesgue integrable, then the integral average of $f$ over $\langle a,b\rangle$ is defined by 
\Eq{IA}{
  \A(f,\langle a,b\rangle):=\int_0^1 f(ta+(1-t)b)dt.
}
One can easily see that the following equality holds:
\Eq{IA+}{
  \langle a,b\rangle
  :=\begin{cases}
   [a,b] &\mbox{ if } a<b, \\[1mm]
   \{a\} &\mbox{ if } a=b, \\[1mm]
   [b,a] &\mbox{ if } a>b,
   \end{cases}
    \qquad\mbox{and}\qquad
   \A(f,\langle a,b\rangle)
  =\begin{cases}
   \dfrac{1}{b-a}\displaystyle\int_a^b f &\mbox{ if } a<b, \\[4mm]
   f(a) &\mbox{ if } a=b, \\[2mm]
   \dfrac{1}{a-b}\displaystyle\int_b^a f &\mbox{ if } a>b.
   \end{cases}
}

The inequalities stated in the following results summarize both Hermite--Hadamard and Ostrow\-ski-type inequalities for the $\Phi$-monotone as well as for the $\Phi$-H\"older settings.

\Thm{GHH}{Let $\Phi\in\E(I)$ and $f:I\to\R$ be locally Lebesgue integrable functions. If $f$ is $\Phi$-monotone, then, for all $u,v,w,z\in I$ with $u\leq w$ and $v\leq z$,
\Eq{GHH}{
  \A(f,\langle u,v\rangle)
  \leq \A(f,\langle w,z\rangle)+\A(\Phi,\langle w-u,z-v\rangle).
}
If $f$ is $\Phi$-H\"older, then, for all $u,v,w,z\in I$,
\Eq{GHH2}{
  \big|\A(f,\langle u,v\rangle)-\A(f,\langle w,z\rangle)\big|
  \leq\A(\Phi\circ|\cdot|,\langle w-u,z-v\rangle).
}}

\begin{proof} Assume first that $f$ is $\Phi$-monotone. Then the inequalities $u\leq w$ and $v\leq z$ imply $tu+(1-t)v\leq tw+(1-t)z$ for $t\in[0,1]$, hence
\Eq{*}{
   f(tu+(1-t)v)\leq f(tw+(1-t)z)+\Phi(t(w-u)+(1-t)(z-v)).
}
Integrating with respect to $t$ over $[0,1]$, we get that \eq{GHH} holds.

In the case when $f$ is $\Phi$-H\"older, for $t\in[0,1]$, we get
\Eq{*}{
   f(tu+(1-t)v)-f(tw+(1-t)z)\leq \Phi(|t(w-u)+(1-t)(z-v)|),\\
   f(tw+(1-t)z)-f(tu+(1-t)v)\leq \Phi(|t(w-u)+(1-t)(z-v)|).
}
Integrating both inequalities with respect to $t$ over $[0,1]$, we get that \eq{GHH2} holds.
\end{proof}

Assuming $\Phi$-monotonicity, we deduce a monotonicity type integral inequality which we will call the \emph{lower and upper Hermite--Hadamard inequalities for $\Phi$-monotone functions}.

\Thm{HH}{Let $\Phi\in\E(I)$ and $f:I\to\R$ be a $\Phi$-monotone. Assume that both functions are locally Lebesgue integrable. Then, for every $x<y$ in $I$, the following two inequalities hold:
\Eq{HH}{
  f(x)-\frac{1}{y-x}\int_0^{y-x} \Phi
  \leq \frac{1}{y-x}\int_x^y f
  \leq f(y)+\frac{1}{y-x}\int_0^{y-x} \Phi.
}
Furthermore, if $\Phi$ is subadditive and nondecreasing, then, for all $x<y$ in $I$,
\Eq{HH1}{
  \sup_{f\in \M_\Phi(I)} \bigg(f(x)-\frac{1}{y-x}\int_x^y f\bigg)
  =\sup_{f\in \M_\Phi(I)} \bigg(\frac{1}{y-x}\int_x^y f-f(y)\bigg)
  =\frac{1}{y-x}\int_0^{y-x} \Phi.
}}

\begin{proof}
Let $x<y$. Then the left and right hand side inequalities in \eq{HH} follow from \eq{GHH} and \eq{IA+} by taking the particular cases $u=v=w=x$, $z=y$ and $u=x$, $v=w=z=y$, respectively. 

Now assume that $\Phi$ is subadditive and nondecreasing. In view of \eq{HH}, it is clear that both supremums are not bigger than the right hand side of \eq{HH1}. To prove the equalities in \eq{HH1}, let $x<y$ in $I$. Then let $\underline{f}$ and $\overline{f}$ be defined by
\Eq{*}{
  \underline{f}(u):=\begin{cases}
         0 &\mbox{if } u\leq x,\\%[2mm]
         -\Phi(u-x) &\mbox{if } x<u,
        \end{cases}
    \qquad\mbox{and}\qquad
 \overline{f}(u):=\begin{cases}
         \Phi(y-u) &\mbox{if } u<y,\\%[2mm]
         0 &\mbox{if } y\leq u.
        \end{cases}
}
By \prp{1}, we can obtain that both $\underline{f}$ and $\overline{f}$ are  $\Phi$-monotone. On the other hand, we have
\Eq{*}{
  \sup_{f\in \M_\Phi(I)} \bigg(f(x)-\frac{1}{y-x}\int_x^y f\bigg)
  \geq \underline{f}(x)-\frac{1}{y-x}\int_x^y \underline{f}
  =\frac{1}{y-x}\int_x^y \Phi(u-x)du=\frac{1}{y-x}\int_0^{y-x} \Phi
}
and 
\Eq{*}{
  \sup_{f\in \M_\Phi(I)} \bigg(\frac{1}{y-x}\int_x^y f-f(y)\bigg)
  \geq \frac{1}{y-x}\int_x^y \overline{f}-\overline{f}(y)
  = \frac{1}{y-x}\int_x^y \Phi(y-u)du=\frac{1}{y-x}\int_0^{y-x} \Phi.
}
These two inequalities together with their reverses imply  that  the equality \eq{HH1} holds.
\end{proof}

\Cor{HH}{Let $p\in[0,1]$, $c\in[0,\infty[\,$ and $f:I\to\R$ be a $c(\cdot)^p$-monotone locally Lebesgue integrable function. Then, for every $x<y$ in $I$, the following two inequalities hold:
\Eq{*}{
  f(x)-\frac{c}{p+1}(y-x)^p
  \leq \frac{1}{y-x}\int_x^y f
  \leq f(y)+\frac{c}{p+1}(y-x)^p .
}
Furthermore, for all $x<y$ in $I$,
\Eq{*}{
  \sup_{f\in \M_\Phi(I)} \bigg(f(x)-\frac{1}{y-x}\int_x^y f\bigg)
  =\sup_{f\in \M_\Phi(I)} \bigg(\frac{1}{y-x}\int_x^y f-f(y)\bigg)
  =\frac{c}{p+1}(y-x)^p.
}}

\begin{proof} Apply the previous statement for the error function $\Phi\in\E(R_+)$ given by $\Phi(t):=ct^p$. This error function is subadditive and nondecreasing, therefore, the second part of the theorem can also be applied.
\end{proof}

\Lem{Phi}{Let $\Psi\in\E(I)$ and assume that the map $t\mapsto\Psi(t)/t$ is locally integrable on $[0,\ell(I)[\,$ and define $\Phi\in\E(I)$ by
\Eq{Phi}{
   \Phi(u)=\Psi(u)+\int_0^u\frac{\Psi(t)}{t}dt \qquad (u\in\,]0,\ell(I)[).
}
Then $\Phi$ is locally integrable and satisfies the following equation:
\Eq{CPhi}{
  \Psi(u) + \frac1u\int_0^u\Phi=\Phi(u) \qquad (u\in\,]0,\ell(I)[).
}}

\begin{proof} By the assumption, we also have that $\Psi$ is locally integrable and hence $\Phi$ defined by \eq{Phi} is also locally integrable. Let $u\in\,]0,\ell(I)[$ be fixed. Then, by Fubini's theorem, we obtain
\Eq{*}{
  \Psi(u) + \frac1u\int_0^u\Phi(s)ds
  &=\Psi(u) + \frac1u\int_0^u \Big(\Psi(s)+\int_0^s\frac{\Psi(t)}{t}dt\Big)ds \\
  &=\Psi(u) + \frac1u\int_0^u \Psi(s)ds
    +\frac1u\int_0^u\int_0^s\frac{\Psi(t)}{t}dtds \\
  &=\Psi(u) + \frac1u\int_0^u \Psi(s)ds
    +\frac1u\int_0^u\int_t^u\frac{\Psi(t)}{t}dsdt \\
  &=\Psi(u) + \frac1u\int_0^u \Psi(s)ds
    +\frac1u\int_0^u(u-t)\frac{\Psi(t)}{t}dt \\
  &=\Psi(u) + \frac1u\int_0^u \Psi(s)ds
    +\int_0^u\frac{\Psi(t)}{t}dt-\frac1u\int_0^u\Psi(t)dt \\
  &=\Psi(u) +\int_0^u\frac{\Psi(t)}{t}dt=\Phi(u).
}
This completes the proof of \eq{CPhi}.
\end{proof}

\Thm{HHLI}{Let $\Psi\in\E(I)$ and assume that the map $t\mapsto\Psi(t)/t$ is locally integrable on $[0,\ell(I)[\,$ and define $\Phi\in\E(I)$ by \eq{Phi}. If $f:I\to\R$ is an upper semicontinuous solution of 
\Eq{HHL}{
  f(u) \leq \frac{1}{v-u}\int_u^v f + \Psi(v-u) \qquad(u,v\in I,\,u<v),
}
then $f$ is $\Phi$-monotone on $I$.}

\begin{proof} Let $f:I\to \R$ be an upper semicontinuous solution of \eq{HHL}. Then $f$ is Lebesgue measurable. Therefore $f$ is upper bounded on any compact subinterval $[u,v]$ of $I$ and hence the Lebesgue integral of $f$ exists and it cannot be $+\infty$. 
On the other hand, \eq{HHL} shows that $(v-u)(f(u)-\Psi(v-u))$ is a lower bound for the 
integral of $f$, thus $f$ has a finite Lebesgue integral over $[u,v]$. 

To prove that $f$ is $\Phi$-monotone, let $x,y\in I$ be fixed with $x<y$. 
Let $USC([x,y])$ denote the family of upper semicontinuous Lebesgue integrable functions and, for $g\in USC([x,y])$, define  
\Eq{*}{
  (Tg)(u):=
  \begin{cases}
  \displaystyle\frac{1}{y-u}\int_u^y g + \Psi(y-u) &\mbox{ if } u\in[x,y[\,,\\[3mm]
  g(y) &\mbox{ if } u=y.
  \end{cases}
}
Then $T$ is a monotone and affine operator which maps $USC([x,y])$ into itself.

Observe that inequality \eq{HHL} for $v=y$ gives that $f\leq Tf$ holds on $[x,y]$. By the monotonicity of $T$, it follows that $Tf\leq T(Tf)=T^2f$. By induction, this yields that $T^{n-1}f\leq T^{n}f$ for all $n\in\N$. Therefore, the sequence $T^nf$ is increasing.

Define $F:[x,y]\to\R$ by 
\Eq{*}{
  F(u):=\sup_{[u,y]} f.
}
By the upper semicontinuity of $f$, we have that $F$ is finite. 

Let $v\in\,]x,y[\,$ be arbitrarily fixed. We prove by induction, for all $n\in\N\cup\{0\}$, that
\Eq{Tn}{
  T^nf(u)\leq
  \begin{cases}
  \Big(\dfrac{v-x}{y-x}\Big)^n (F(x)-F(v)) + F(v) +\Phi(y-u)
     &\mbox{if } u\in[x,v[\,,\\[2mm]
  F(v) + \Phi(y-u) &\mbox{if } u\in[v,y]\,.
  \end{cases}
}
If $n=0$, then $T^nf(u)=f(u)$ and \eq{Tn} simplifies to
\Eq{*}{
  f(u)\leq
  \begin{cases}
  F(x) + \Phi(y-u)  &\mbox{if } u\in[x,v[\,,\\[2mm]
  F(v) + \Phi(y-u) &\mbox{if } u\in[v,y]\,,
  \end{cases}
}
which follows from the definition of $F$ and the nonnegativity of $\Phi$.

Now assume that \eq{Tn} holds for some $n$. For $u\in[x,v[\,$, using the first and second inequalities in \eq{Tn} on the intervals $[u,v]$ and $[v,y]$, respectively, and finally the assertion of \lem{Phi}, we get 
\Eq{*}{
  (T^{n+1}f)(u)
  &=\frac{1}{y-u}\int_u^y T^nf + \Psi(y-u)
  =\frac{1}{y-u}\Big(\int_u^v T^nf + \int_v^y T^nf\Big) + \Psi(y-u)\\
  &\leq\frac{v-u}{y-u} \bigg(\Big(\dfrac{v-x}{y-x}\Big)^n (F(x)-F(v)) + F(v)\bigg) +\frac{1}{y-u}\int_u^v\Phi(y-t)dt \\
  &\qquad + \frac{y-v}{y-u} F(v) + \frac{1}{y-u}\int_v^y\Phi(y-t)dt+\Psi(y-u)\\
  &=\frac{v-u}{y-u}\Big(\dfrac{v-x}{y-x}\Big)^n (F(x)-F(v)) + F(v) + \frac{1}{y-u}\int_u^y\Phi(y-t)dt+\Psi(y-u)\\
  &\leq\bigg(\frac{v-x}{y-x}\bigg)^{n+1} (F(x)-F(v)) + F(v) + \frac{1}{y-u}\int_0^{y-u}\Phi+\Psi(y-u)\\
  &=\bigg(\frac{v-x}{y-x}\bigg)^{n+1} (F(x)-F(v)) + F(v) + \Phi(y-u).
}
If $u\in[v,y[\,$, then $[u,y]\subseteq[v,y]$ and hence we shall only need the second inequality for $T^nf$ on $[u,y]$. 
Using this estimate and the assertion of \lem{Phi}, we get
\Eq{*}{
  (T^{n+1}f)(u)
  &=\frac{1}{y-u}\int_u^y T^nf + \Psi(y-u)
  \leq F(v) + \frac{1}{y-u}\int_u^y\Phi(y-t)dt+\Psi(y-u)\\
  &= F(v) + \frac{1}{y-u}\int_0^{y-u}\Phi+\Psi(y-u)
  = F(v) + \Phi(y-u).
}
and the trivial inequality $(T^{n+1}f)(y)=f(y)\leq F(v) + \Phi(y-u)$ shows  the desired inequality is true for $u=y$. 

In view of the monotonicity of the sequence $T^nf$, we have that $f\leq T^nf$ for all $n\in\N$. Thus, \eq{Tn} implies, for all $n\in\N$, that
\Eq{*}{
  f(u)\leq
  \begin{cases}
  \Big(\dfrac{v-x}{y-x}\Big)^n (F(x)-F(v)) + F(v) +\Phi(y-u)
     &\mbox{if } u\in[x,v[\,,\\[2mm]
  F(v) + \Phi(y-u) &\mbox{if } u\in[v,y]\,.
  \end{cases}
}
Upon taking the limit $n\to\infty$, the above inequality yields, for all $v\in[x,y[\,$, that
\Eq{*}{
  f(u)\leq  F(v) + \Phi(y-u)=\sup_{[v,y]}f + \Phi(y-u) \qquad u\in[x,y].
}
Now, taking the limit $v\to y$, the upper semicontinuity of $f$ yields, that
\Eq{*}{
  f(u)\leq  f(y) + \Phi(y-u) \qquad u\in[x,y].
}
In particular, this inequality holds for $u=x$, which completes the proof of the $\Phi$-monotonicity of $f$.
\end{proof}

The following result is a counterpart of \thm{HHLI}. It can be proved directly in an analogous way, however, we will deduce it from this theorem using a sign transformation.

\Thm{HHRI}{Let $\Psi\in\E(I)$ and assume that the map $t\mapsto\Psi(t)/t$ is locally integrable on $[0,\ell(I)[\,$ and define $\Phi\in\E(I)$ by \eq{Phi}. If $f:I\to\R$ is a lower semicontinuous solution of 
\Eq{HHR}{
  \frac{1}{v-u}\int_u^v f  \leq f(v)+ \Psi(v-u) \qquad(u,v\in I,\,u<v),
}
then $f$ is $\Phi$-monotone on $I$.}

\begin{proof} Assume that $f:I\to\R$ satisfies \eq{HHR}. Define $g(x):=-f(-x)$ for $x\in J:=-I$. Then $\ell(I)=\ell(J)$ hence $\Phi,\Psi\in\E(J)$ and $g$ is lower semicontinuous over $J$. On the other hand, for $u,v\in J$ with $u<v$, we have that $-u,-v\in I$ with $-v<-u$. Applying \eq{HHL} for these variables, we
obtain
\Eq{*}{
  \frac{1}{-u-(-v)}\int_{-v}^{-u} f  \leq f(-u)+ \Psi(-u-(-v)) 
  \qquad(u,v\in J,\,u<v).
}
This, after replacing $f(x)$ by $-g(-x)$, gives 
\Eq{*}{
  -\frac{1}{v-u}\int_{u}^{v}g  \leq -g(u)+ \Psi(v-u) 
  \qquad(u,v\in J,\,u<v).
}
From here, we can see that $g$ satisfies the inequality \eq{HHL} (wherein $f$ is replaced by $g$). According to \thm{HHLI}, we can conclude that $g$ is $\Phi$-monotone on $J$. This, with an analogous substitutions, implies that $f$ is also $\Phi$-monotone.
\end{proof}

\Cor{HHLI}{Let $p\in\,]0,1]$, $c\in[0,\infty[\,$. If $f:I\to\R$ is an upper semicontinuous solution of 
\Eq{*}{
  f(u) \leq \frac{1}{v-u}\int_u^v f + c(v-u)^p \qquad(u,v\in I,\,u<v),
}
then $f$ is $\frac{c(p+1)}{p}(\cdot)^p$-monotone on $I$. In particular, $f$ is increasing if $p>1$.}

\begin{proof} By our assumption, $f$ satisfies \eq{HHL} with $\Psi$ defined by $\Psi(u):=cu^p$. In order to apply the previous theorem, we have to compute $\Phi$ which is given by \eq{Phi}.
\Eq{*}{
  \Phi(u)=\Psi(u)+\int_0^u\frac{\Psi(t)}{t}dt
  =cu^p + \int_0^uct^{p-1}dt=cu^p + \frac{c}{p}u^p=\frac{c(p+1)}{p}u^p.
}
Thus, by \thm{HHLI}, $f$ is $\frac{c(p+1)}{p}(\cdot)^p$-monotone on $I$.

In the case $p>1$, the $\frac{c(p+1)}{p}(\cdot)^p$-monotonicity implies that $f$ is in fact increasing.
\end{proof}

The next assertion is a counterpart of \cor{HHLI}. Its proof can be obtained from \thm{HHRI} exactly in the same way as the previous corollary was deduced from \thm{HHLI}.

\Cor{HHRI}{Let $p\in\,]0,1]$, $c\in[0,\infty[\,$. If $f:I\to\R$ is a lower semicontinuous solution of 
\Eq{*}{
  \frac{1}{v-u}\int_u^v f\leq f(v) + c(v-u)^p \qquad(u,v\in I,\,u<v),
}
then $f$ is $\frac{c(p+1)}{p}(\cdot)^p$-monotone on $I$. In particular, $f$ is increasing if $p>1$.}

\section{Ostrowski- and Hermite--Hadamard-type inequalities for $\Phi$-H\"older functions}

The result stated in the next theorem will be called an Ostrowski-type inequality  for $\Phi$-H\"older functions.

\Thm{OI}{Let $\Phi\in\E(I)$ and $f:I\to\R$ be a $\Phi$-H\"older. Assume that both functions are locally Lebesgue integrable. Then, for every $x<y$ in $I$, the following two inequality hold:
\Eq{OI}{
  \bigg|f(p)-\frac{1}{y-x}\int_x^y f\bigg| 
  \leq \frac{1}{y-x}\bigg(\int_0^{p-x}\Phi+\int_0^{y-p}\Phi\bigg)
  \qquad(p\in[x,y]).
}
Furthermore, if $\Phi$ is subadditive and nondecreasing with $\Phi(0)=0$, then, for all $x<y$ in $I$,
\Eq{OI1}{
  \sup_{f\in \H_\Phi(I)} \bigg|f(p)-\frac{1}{y-x}\int_x^y f\bigg|
  =\frac{1}{y-x}\bigg(\int_0^{p-x}\Phi+\int_0^{y-p}\Phi\bigg)\qquad(u\in[x,y]).
}}

\begin{proof} Let $p\in[x,y]$ be fixed. Applying the second assertion of \thm{GHH} in the particular case $u=v=p$, $w=x$, and $z=y$, the inequality \eq{GHH2} yields
\Eq{*}{
  \big|\A(f,\{p\})-\A(f,[x,y])\big|
  &\leq\A(\Phi\circ|\cdot|,[x-p,y-p])\\
  &=\frac{1}{y-x}\int_{x-p}^{y-p} \Phi(|t|) dt
  =\frac{1}{y-x}\bigg(\int_0^{p-x}\Phi+\int_0^{y-p}\Phi\bigg).
}
which reduces to \eq{OI}.

Assume now that $\Phi$ is subadditive and nondecreasing and let $x<y$ be fixed. For $p\in[x,y]$, define the function $\Phi_p:I\to\R$ by \eq{Pp}. Then, by \prp{HH}, we have that $\Phi_p$ is a $\Phi$-H\"older function. Therefore, 
\Eq{*}{
  \sup_{f\in \H_\Phi(I)} \bigg|f(p)-\frac{1}{y-x}\int_x^y f\bigg|
  &\geq \bigg|\Phi_p(p)-\frac{1}{y-x}\int_x^y \Phi_p\bigg|\\
  &=\frac{1}{y-x}\int_x^y \Phi_p
  =\frac{1}{y-x}\bigg(\int_0^{p-x}\Phi+\int_0^{y-p}\Phi\bigg).
}
This inequality together with its reverse imply that the equality in \eq{OI1} is valid.
\end{proof}

The inequalities obtained in the particular cases $p=x$ and $p=y$ will be called the \emph{lower and upper Hermite--Hadamard-type inequalities for $\Phi$-H\"older functions}. 

\Cor{HHH}{Let $q\in[0,1]$, $c\in[0,\infty[\,$ and $f:I\to\R$ be a $c(\cdot)^q$-H\"older locally Lebesgue integrable function. Then, for every $x<y$ in $I$,
\Eq{*}{
  \bigg|f(p)-\frac{1}{y-x}\int_x^y f\bigg| 
  \leq \frac{c}{q+1}((p-x)^q+(y-p)^q) \qquad(p\in[x,y]).
}
Furthermore, for all $x<y$ in $I$,
\Eq{*}{
  \sup_{f\in \H_\Phi(I)} \bigg|f(p)-\frac{1}{y-x}\int_x^y f\bigg|
  =\frac{c}{q+1}((p-x)^q+(y-p)^q) \qquad(p\in[x,y]).
}}

\begin{proof} Apply the previous statement for the error function $\Phi\in\E(R_+)$ given by $\Phi(t):=ct^p$.
\end{proof}

\Thm{HHHLI}{Let $\Psi\in\E(I)$ and assume that the map $t\mapsto\Psi(t)/t$ is locally integrable on $[0,\ell(I)[\,$ and define $\Phi\in\E(I)$ by \eq{Phi}. If $f:I\to\R$ is a continuous solution of 
\Eq{HHHL}{
  \bigg|f(u)-\frac{1}{v-u}\int_u^v f\bigg| \leq \Psi(v-u) \qquad(u,v\in I,\,u<v),
}
then $f$ is $\Phi$-H\"older on $I$.}

\begin{proof} If $f$ satisfies \eq{HHHL}, then $f$ and $-f$ fulfil \eq{HHL}. Therefore, in view of \thm{HHLI}, we obtain that $f$  and $-f$ are $\Phi$-monotone, which implies that $f$ is $\Phi$-H\"older on $I$.
\end{proof}

The proofs of the following results are similar to that of the previous theorem, they are left to the reader.

\Thm{HHHRI}{Let $\Psi\in\E(I)$ and assume that the map $t\mapsto\Psi(t)/t$ is locally integrable on $[0,\ell(I)[\,$ and define $\Phi\in\E(I)$ by \eq{Phi}. If $f:I\to\R$ is a continuous solution of 
\Eq{*}{
  \bigg|f(v)-\frac{1}{v-u}\int_u^v f\bigg| \leq \Psi(v-u) \qquad(u,v\in I,\,u<v),
}
then $f$ is $\Phi$-H\"older on $I$.}

\Cor{HHHLI}{Let $p\in\,]0,1]$, $c\in[0,\infty[\,$. If $f:I\to\R$ is a continuous solution of 
\Eq{*}{
  \bigg|f(u)-\frac{1}{v-u}\int_u^v f\bigg| \leq c(v-u)^p \qquad(u,v\in I,\,u<v),
}
then $f$ is $\frac{c(p+1)}{p}(\cdot)^p$-H\"older on $I$. In particular, $f$ is constant if $p>1$.}

\Cor{HHHRI}{Let $p\in\,]0,1]$, $c\in[0,\infty[\,$. If $f:I\to\R$ is a continuous solution of 
\Eq{*}{
  \bigg|f(v)-\frac{1}{v-u}\int_u^v f\bigg| \leq c(v-u)^p \qquad(u,v\in I,\,u<v),
}
then $f$ is $\frac{c(p+1)}{p}(\cdot)^p$-H\"older on $I$. In particular, $f$ is constant if $p>1$.}

\section{A concluding remark}

The superadditivity of nonpositive functions can easily be characterized
according to the next statement.

\Prp{pppp}{Let $I\subseteq\R_+$ with $\inf I=0$. Let $f: I\to \R$ be a nonpositive function. Then $f$ is superadditive if and only if $f$ is $(-f)$-monotone.}

\begin{proof}
Assume first that $f$ is $(-f)$-monotone. Let $x,y\in I$ with $x+y\in I$. Then the $(-f)$-monotonicity of $f$ implies
\Eq{*}
{f(x)\leq f(x+y)+(-f)(y).
}
which yields that $f$ is superadditive.

Conversely, assume that $f$ is superadditive. Then, for any $x,y\in I$ with $x\leq y$,
\Eq{*}
{f(y)=f(x+(y-x))\geq f(x)+f(y-x)=f(x)-(-f)(y-x).
}
Which establishes the $(-f)$-monotonicity of $f$.
\end{proof}

%\subsection*{Acknowledgement} The authors are grateful for the detailed report of the anonymous referee.

%\bibliography{publ,funcequ}
%\bibliographystyle{plain}

\end{document}